\documentclass[11pt]{article}
\usepackage{amsmath}
\usepackage{amssymb}
\usepackage{amsfonts}
\usepackage{color}
\allowdisplaybreaks

\newtheorem{Theorem}{{\bf Theorem}}[section]
\newtheorem{Lemma}{{\bf Lemma}}[section]
\newtheorem{Proposition}{{\bf Proposition}}[section]
\newtheorem{Corollary}{{\bf Corollary}}[section]
\newtheorem{Remark}{{\bf Remark}}[section]
\newtheorem{Example}{{\bf Example}}[section]
\newtheorem{Definition}{{\bf Definition}}[section]

\newenvironment{theorem}{\begin{Theorem}$\!\!\!$}{\end{Theorem}}
\newenvironment{lemma}{\begin{Lemma}$\!\!\!$}{\end{Lemma}}
\newenvironment{proposition}{\begin{Proposition}$\!\!\!$}{\end{Proposition}}
\newenvironment{corollary}{\begin{Corollary}$\!\!\!$}{\end{Corollary}}
\newenvironment{remark}{\begin{Remark}$\!\!\!$}{\end{Remark}}
\newenvironment{example}{\begin{Example}$\!\!\!$}{\end{Example}}

\numberwithin{equation}{section}


\begin{document}
\title{Moment conditions and lower bounds in expanding solutions of \\
wave equations with double damping terms}

\author{Ryo Ikehata\thanks{ikehatar@hiroshima-u.ac.jp}\\
{\small Department of Mathematics, Graduate School of Education, Hiroshima University} \\
{\small Higashi-Hiroshima 739-8524, Japan}
\and Hironori Michihisa\thanks{Corresponding author: hi.michihisa@gmail.com}\\ 
{\small Department of Mathematics, Graduate School of Science, Hiroshima University} \\
{\small Higashi-Hiroshima 739-8526, Japan}
}
\date{}

\maketitle

\begin{abstract}
In this report we obtain higher order asymptotic expansions of solutions to 
wave equations with frictional and viscoelastic damping terms. 
Although the diffusion phenomena are dominant, 
differences between the solutions we deal with and those of heat equations can be seen by comparing the second order expansions of them. 
In order to analyze such effects we consider the weighted $L^1$ initial data. 
We also give some lower bounds which show the optimality of obtained expansions. 
\end{abstract}

\footnote[0]{\hspace{-2em} 2010 {\it Mathematics Subject Classification}. primary 35B40, 35C20; secondary 35L05, 35K08}
\footnote[0]{\hspace{-2em} {\it Keywords and Phrases}: 
Wave equation; 
moment condition; 
asymptotic expansion; 
lower bounds estimates; 
diffusion phenomena; 
double damping terms
}

\section{Introduction}
We consider the $n$-dimensional Cauchy problem of wave equations with frictional and viscoelastic damping terms 
\begin{align}
\label{1.1}
\begin{cases}
u_{tt}-\Delta u+u_t-\Delta u_t=0, 
& (t,x)\in (0,\infty)\times{\bf R}^n, \\
u(0,x)=u_0(x),\quad u_t(0,x)=u_1(x), 
& x\in{\bf R}^n, 
\end{cases}
\end{align}
where $n\ge1$ and $u_0$, $u_1\in L^2({\bf R}^n)\cap L^{1,\gamma}({\bf R}^n)$ with $\gamma\ge0$. 
Here the weighted $L^1$-space $L^{1,\gamma}({\bf R}^n)$ is defined by 
\[
L^{1,\gamma}({\bf R}^n)
:=\left\{f\in L^1({\bf R}^n)\biggr| 
\|f\|_{1,\gamma}:=\int_{{\bf R}^n}(1+|x|)^k |f(x)|\,dx<\infty \right\}.
\]

Before referring to some previous papers (e.g.,\cite{IS} and \cite{D}) dealing directly with the equation~\eqref{1.1}, 
we mention several results on classical wave equations with damping terms. 

Many mathematicians have studied the asymptotic profiles and decay estimates of solutions to wave equations with damping terms. 
For example, we can refer the reader to 
\cite{H}, \cite{HO}, \cite{K}, \cite{MN}, \cite{M}, \cite{Mi}, \cite{Na}, \cite{Ni}, \cite{SW}, \cite{T} 
on damped wave equations, and \cite{DR}, \cite{I}, \cite{IO}, \cite{IT}, \cite{ITY}, \cite{Mi2}, \cite{P}, \cite{S}
on strongly damed wave equations. 

To the best of authors' knowledge, 
the diffusive structure strongly appears in wave equations with the damping term $u_t$. 
It is well-known that the solution $u$ to the following damped wave equation 
\begin{align*}
\begin{cases}
u_{tt}-\Delta u +u_t=0, & t>0, \quad x\in{\bf R}^n, \\
u(0,x)=u_0(x), \quad u_t(0,x)=u_1(x), & x\in{\bf R}^n,
\end{cases}
\end{align*}
behaves like a solution $v$ to the corresponding heat equation 
\begin{align*}
\begin{cases}
v_t-\Delta v=0, & t>0, \quad x\in{\bf R}^n, \\
v(0,x)=u_0(x)+u_1(x), & x\in{\bf R}^n,
\end{cases}
\end{align*}
as time goes to infinity. 
While, the solution of the strongly damped wave equation
\[
u_{tt}-\Delta u - \Delta u_t = 0
\]
behaves like the so-called diffusion wave as $t \to \infty$. Recently, this diffusion wave property is successively studied (for example) in \cite{I}, \cite{IT}, \cite{ITY}, \cite{Mi} and \cite{S}. 

On the other hand, quite recently Ikehata-Sawada \cite{IS} consider the equation \eqref{1.1}, and investigated which term ($u_t$ or $-\Delta u_t$) gives stronger effects on the asymptotic profile of the solution. They concluded that the term $u_t$ is more dominant as $t \to \infty$. Indeed, by considering the Fourier transform of $u_t$ and $-\Delta u_t$ with respect to the spatial variable, this observation seems natural. 
The diffusive structure is essentially seen in the low frequency region $\{|\xi|\ll1\}$ in the Fourier space, and so the damping term $\hat{u}_t$ is stronger than the viscoelastic one $|\xi|^2 \hat{u}_t$. 
In \cite{IS}, they obtained the following asymptotic estimates for the solution to \eqref{1.1} with $[u_0,u_1]\in (H^1({\bf R}^n)\cap L^{1,1}({\bf R}^n))\times(L^2({\bf R}^n)\cap L^{1,1}({\bf R}^n))$: 
\[
\biggr\|u(t)-(P_{00}+P_{01})G(t)\biggr\|_2
\le Ct^{-\frac{n}{4}-\frac{1}{2}}
\big(\|u_0\|_{1,1}+\|u_1\|_{1,1}+\|u_0\|_{H^2}+\|u_1\|_2\big), 
\qquad
t\gg1.
\]
Here 
\[
P_{0j}
:=\int_{{\bf R}^n}
u_j(x)\,dx, 
\qquad
j=0,1, 
\]
and the function $G(t,x)$ is the Gauss kernel 
\[
G(t,x):=(4\pi t)^{-\frac{n}{4}}\exp\left(-\frac{|x|^2}{4t}\right).
\]
They have additionally obtained the higher order (up to the first order) asymptotic profiles of the solution to \eqref{1.1} after identifying the above leading term (the zero-th order expansion) only in the one  dimensional case. However, since the solution to \eqref{1.1} is very close to the one of the heat equation as $t \to \infty$, 
one can expect further asymptotic expansions of the solution to \eqref{1.1} under more heavy moment condition on the initial data in all space dimensional cases (see e.g., \cite{IKM}). 
When we want to observe the difference between them, 
we have to get the asymptotic expansions of the solution to \eqref{1.1} higher than the second order. 
In this connection, after \cite{IS}, D'Abbicco \cite{D} derives the several decay estimates of the solution to problem \eqref{1.1}, 
and applied them to nonlinear problems with three types of nonlinearities like $|u|^{p}$, $|u_{t}|^{p}$ and $|\nabla u|^{p}$ in order to investigate the so-called critical exponent for $p$. 
Recently, independently from \cite{D}, 
Ikehata-Takeda \cite{IT-1, IT-2} have caught the asymptotic profile of solutions and critical exponent of the power $p$ of the nonlinearity $|u|^{p}$ to the semilinear problem of \eqref{1.1} in the low dimensional case $n = 1,2,3$. 
Concerning the higher order asymptotic expansions of solutions to problem \eqref{1.1} it seems that we still do not have any previous research manuscripts.   

The purpose of this paper is to report higher order asymptotic expansions and 
some lower bounds of solutions to problem \eqref{1.1}, which imply their optimality. \\

This paper is organized as follows. 
In section 2, we collect several preliminary results and notation, which will be used throughout this paper.  
In section 3, we introduce two important functions and their properties, which determine the desired asymptotic profiles of the solution to problem \eqref{1.1}. Main results are stated in section 4, and the proofs of the results will be given in section 5. 
In section 6, we shall reconsider the heat equation to compare with the profiles of \eqref{1.1}. Appendix will be given in section 7 .

\section{Notation}
In this section, we will introduce some notation and preliminary knowledges, which will be used throughout this paper.

We write ${\bf N}$ as the set of all positive integers, and put  ${\bf N}_0:={\bf N}\cup\{0\}$. For $r\in{\bf R}$, 
we denote the maximal integer $N$ satisfying $N\le r$ by $[r]$. When the Fourier transform of $f$ can be defined, 
it is defined by 
\[
\mathcal{F}[f](\xi)
=\hat{f}(\xi)
:=\int_{{\bf R}^n} e^{-ix\cdot\xi}f(x)\,dx.
\] 
Note that the Fourier transform of the Gauss kernel is 
\[
\mathcal{F}[G(t,\cdot)](\xi)
=e^{-t|\xi|^2}
\]
under this definition. 

Applying the Fourier transform to \eqref{1.1}, 
we see that
\begin{align*}
\begin{cases}
\hat{u}_{tt}+|\xi|^2 \hat{u}+(1+|\xi|^2)\hat{u}_t=0, 
&  (t,\xi)\in (0,\infty)\times{\bf R}^n, \\
\hat{u}(0,\xi)=\hat{u}_0(\xi),\qquad \hat{u}_t(0,\xi)=\hat{u}_1(\xi),  
& \xi\in{\bf R}^n.
\end{cases}
\end{align*}
Solving the corresponding characteristic equation: 
\[
\lambda^{2} + (1+|\xi|^2)\lambda+|\xi|^2=0,
\]
we have 
\[
\lambda_\pm=\frac{-(1+|\xi|^2)\pm\big|1-|\xi|^2\big|}{2},
\]
that is, 
\begin{align*}
\begin{cases}
\lambda_+=-|\xi|^2,
\quad
\lambda_-=-1 
& \mbox{if}
\qquad
|\xi|\le1, \\
\lambda_+=-1,
\quad
\lambda_-=-|\xi|^2 
& \mbox{if}
\qquad
|\xi|\ge1.
\end{cases}
\end{align*}
Thus we have 
\begin{align}
\hat{u}(t,\xi)
& =\frac{1}{1-|\xi|^2}e^{-t|\xi|^2}(\widehat{u_0}+\widehat{u_1})
-\frac{1}{1-|\xi|^2}e^{-t}(|\xi|^2 \widehat{u_0}+\widehat{u_1}) \label{2.1}\\
& =\frac{e^{-t|\xi|^2}-|\xi|^2 e^{-t}}{1-|\xi|^2}\widehat{u_0}
+\frac{e^{-t|\xi|^2}-e^{-t}}{1-|\xi|^2}\widehat{u_1} \label{2.2}\\
& =\frac{1}{|\xi|^2-1}e^{-t}(|\xi|^2 \widehat{u_0}+\widehat{u_1})
-\frac{1}{|\xi|^2-1}e^{-t|\xi|^2}(\widehat{u_0}+\widehat{u_1}) \label{2.3}.
\end{align}
In \cite{ITY}, it is shown that there exists a unique weak solution to \eqref{1.1} in the class 
\[
u\in C^1([0,\infty);L^2({\bf R}^n))
\cap 
C([0,\infty);H^1({\bf R}^n))
\]
in the case when 
\[
[u_0,u_1]
\in H^1({\bf R}^n)\times L^2({\bf R}^n). 
\]
The right-hand sides of \eqref{2.1}-\eqref{2.3} are equivalent for each other, and define the $L^2({\bf R}_{\xi}^{n})$-function for each $t\ge0$ 
even if $u_0$, $u_1\in L^2({\bf R}^n)$. 
So we regard the function $\hat{u}(t,\xi)$ defined by \eqref{2.1}-\eqref{2.3} as (the Fourier transform of) the solution to \eqref{1.1} even when $u_0$, $u_1\in L^2({\bf R}^n)$. 
In addition to \eqref{2.1}-\eqref{2.3}, 
we can give another form of the solution to \eqref{1.1}. 
When we read \eqref{1.1} as 
\[
(u_t+u)_t-\Delta(u_t+u)=0, 
\qquad 
(u_t+u)(0,x)=u_0(x)+u_1(x), 
\]
an equivalent representation is obtained:
\begin{align}
\label{2.4}
\hat{u}(t,\xi)
=e^{-t}\widehat{u_0}
+\left(
e^{-t}
\int_0^t e^{-s(|\xi|^2-1)}\,ds
\right)
(\widehat{u_0}+\widehat{u_1}).
\end{align}
Expression~\eqref{2.4} shows that $\{|\xi|=1\}$ is not a singular set and that $u \in C([0,\infty);L^{2}({\bf R}^{n}))$ if $u_0$, $u_1\in L^2({\bf R}^n)$.

Finally, for $f\in L^{1,\gamma}({\bf R}^n)$, 
we set 
\[
M_\alpha (f)
:=\frac{(-1)^{|\alpha|}}{\alpha!}
\left(
\int_{{\bf R}^n} x^\alpha f(x)\,dx
\right), 
\qquad
|\alpha|\le[\gamma],
\]
\noindent
and one also defines 
\[
F^v(\xi):=\frac{1}{1-|\xi|^2}\hat{v}(\xi).
\]
The function $F^{v}(\xi)$ is quite important to get the asymptotic profiles of the solution $\hat{u}(t,\xi)$ to problem \eqref{1.1} (see Lemma~\ref{lem:5.2} below).

\section{Asymptotic profiles}
Let $k\in{\bf N}_0$. 
Here we treat a suitable function $v$ that belongs to some subspaces of $L^1$-space with sufficient weights in order to get higher order expansions of $F^v$. 

We define the following functions:
\begin{align*}
A_k^v(\xi)=
\begin{cases}
\displaystyle{
\sum_{j=0}^{k/2}
\left(
|\xi|^{k-2j}
\sum_{|\alpha|\le2j}
M_\alpha(v)(i\xi)^\alpha
\right)
},
& k\equiv0\pmod{2}, \\[22pt]
\displaystyle{
\sum_{j=0}^{(k-1)/2}
\left(
|\xi|^{k-1-2j}
\sum_{|\alpha|\le2j+1}
M_\alpha(v)(i\xi)^\alpha
\right)
},
& k\equiv1\pmod{2}.
\end{cases} 
\end{align*}
For $k\in{\bf N}_0$, put 
\begin{align*}
B_k^v(\xi):=
\begin{cases}
\displaystyle{
\sum_{j=0}^{k/2}
\left(
|\xi|^{k-2j}
\sum_{|\alpha|=2j}
M_\alpha(v)(i\xi)^\alpha
\right)
},
& k\equiv0\pmod{2},\\[22pt]
\displaystyle{
\sum_{j=0}^{(k-1)/2}
\left(
|\xi|^{k-1-2j}
\sum_{|\alpha|=2j+1}
M_\alpha(v)(i\xi)^\alpha
\right)
},
& k\equiv1\pmod{2}. 
\end{cases}
\end{align*}
Then the functions defined above have the following properties: 
\begin{itemize}
\item[(A)] 
Let $k\in{\bf N}_0$.  
It holds that 
\begin{align*}
A_k^v(\xi)=A_{k-1}^v(\xi)+B_k^v(\xi) 
\end{align*}
for $\xi\in{\bf R}^n$. 
For convenience, here we define $A_{-1}^v(\xi)\equiv0$; 
\item[(B)] 
Let $2\le k\in{\bf N}$. 
It follows that 
\begin{align*}
B_k^v(\xi)
=|\xi|^2 B_{k-2}^v(\xi)
+\sum_{|\alpha|=k}M_\alpha(v)(i\xi)^\alpha
\end{align*}
for $\xi\in{\bf R}^n$;
\item[(C)] 
Let $k\in{\bf N}_0$. 
For any $c>0$ and $\xi\in{\bf R}^n$, 
\begin{align*}
B_k^v\left(\frac{\xi}{c}\right)
=c^{-k} B_k^v(\xi).
\end{align*} 
\end{itemize}

Finally, we can see that  
$A_k^v(\xi)e^{-t|\xi|^2}$ and $B_k^v(\xi)e^{-t|\xi|^2}$ are Fourier transform of 
\begin{align*}
\begin{cases}
\displaystyle{
\sum_{j=0}^{k/2}
\sum_{|\alpha|\le2j}
M_\alpha(v)
(-\Delta)^{(k-2j)/2}
\partial_x^\alpha 
G(t,x)
},
& k\equiv0\pmod{2}, \\[22pt]
\displaystyle{
\sum_{j=0}^{(k-1)/2}
\sum_{|\alpha|\le2j+1}
M_\alpha(v)
(-\Delta)^{(k-1-2j)/2}
\partial_x^\alpha
G(t,x)
},
& k\equiv1\pmod{2},
\end{cases} 
\end{align*}
\begin{align*}
\begin{cases}
\displaystyle{
\sum_{j=0}^{k/2}
\sum_{|\alpha|=2j}
M_\alpha(v)
(-\Delta)^{(k-2j)/2}
\partial_x^\alpha 
G(t,x)
},
& k\equiv0\pmod{2}, \\[22pt]
\displaystyle{
\sum_{j=0}^{(k-1)/2}
\sum_{|\alpha|=2j+1}
M_\alpha(v)
(-\Delta)^{(k-1-2j)/2}
\partial_x^\alpha
G(t,x)
},
& k\equiv1\pmod{2},
\end{cases} 
\end{align*}
respectively.

\section{Main Results}
In this section, we introduce our main results, which show higher order asymptotic expansions of the solution to \eqref{1.1} with weighted $L^1$ initial data. 
\begin{theorem}
\label{thm:1}
Let $n\ge1$ and $\gamma\ge0$, and let $\hat{u}$ be the function defined by \eqref{2.1}-\eqref{2.3} 
with $u_0$, $u_1\in L^2({\bf R}^{n})\cap L^{1,\gamma}({\bf R}^{n})$. 
Then, it holds that 
\begin{align}
\label{4.1}
\biggr\|\hat{u}(t)-A_{[\gamma]}^v e^{-t|\xi|^2}\biggr\|_2 
\le C(1+t)^{-\frac{n}{4}-\frac{\gamma}{2}} \|v\|_{1,\gamma}
+Ce^{-ct}\big(\|u_0\|_2+\|u_1\|_2+\|v\|_{1,[\gamma]}\big)
\end{align}
for $t\ge1$, where $v:=u_0+u_1$. 
Here $C>0$ and $c>0$ are constants independent of $t$, $u_0$ and $u_1$. 
\end{theorem}

\begin{remark}
{\rm We can obtain corresponding results to Theorem~\ref{thm:1} for more general equations 
\[
u_{tt}-\Delta u+\alpha u_t+\beta(-\Delta) u_{t}=0 
\]
by changing $A_k$ according to positive constants $\alpha$ and $\beta$. 
Their proofs may not be easy because solution formulae will be more complicated than \eqref{2.1}-\eqref{2.3}. 
At least \eqref{2.4} does not hold.}
\end{remark}

\begin{theorem}
\label{thm:2}
Let $n\ge1$ and $k\in\bf{N}_0$, and let $\hat{u}$ be the function defined by \eqref{2.1}-\eqref{2.3} 
with $u_0$, $u_1\in L^2({\bf R}^n)\cap L^{1,k}({\bf R}^n)$. 
Then, there exists a constant $\delta\ge1$ depending on $u_0$ and $u_1$ such that 
\begin{equation}
\label{4.2}
\begin{split}
\biggr\|\hat{u}(t)-A_{k-1}^v e^{-t|\xi|^2}\biggr\|_2 
\ge \frac{1}{2} 
\biggr\|B_k^v e^{-|\xi|^2}\biggr\|_{L^2(|\xi|\le1/2)}
t^{-\frac{n}{4}-\frac{k}{2}}
\end{split}
\end{equation}
for $t\ge\delta$, where $v:=u_0+u_1$. 
\end{theorem}

By \eqref{4.1} and \eqref{4.2}, 
we can easily obtain the following corollary since 
\[
\biggr\|B_k^v e^{-t|\xi|^2}\biggr\|_2
\le Ct^{-\frac{n}{4}-\frac{k}{2}},
\qquad
t>0.
\]
Here $C>0$ is a constant depending on the momentum of $v$. 
Inequality~\eqref{4.3} below implies that the obtained expansions are optimal. 
\begin{corollary}
\label{cor:4.1}
Let $n\ge1$ and $k\in\bf{N}_0$, and let $\hat{u}$ be the function defined by \eqref{2.1}-\eqref{2.3} 
with $u_0$, $u_1\in L^2({\bf R}^n)\cap L^{1,k}({\bf R}^n)$. Then there exist constants $C > 0$ and $\delta\ge1$ depending on $u_0$ and $u_1$ such that 
\begin{align}
\label{4.3}
\frac{1}{2} 
\biggr\|B_k^v e^{-|\xi|^2}\biggr\|_{L^2(|\xi|\le1/2)}
t^{-\frac{n}{4}-\frac{k}{2}}
\le \biggr\|\hat{u}(t)-A_{k-1}^v e^{-t|\xi|^2}\biggr\|_2 
\le Ct^{-\frac{n}{4}-\frac{k}{2}}
\end{align}
for $t\ge\delta$, where $v:=u_0+u_1$. 
\end{corollary}

To check the positivity of the quantity $\biggr\|B_k^v e^{-|\xi|^2}\biggr\|_{L^2(|\xi|\le1/2)}$ in Corollary~\ref{cor:4.1}, the following two propositions are useful.
\begin{proposition}
\label{prop:1}
Let $n = 1$ and $v\in L^{1,k}(\bf{R})$ with $k\in\bf{N}_0$. 
Then it holds that 
\begin{align*}
\biggr\|B_k^v e^{-|\xi|^2}\biggr\|_{L^2(|\xi|\le1/2)}=
\begin{cases}
\displaystyle{
\left(
2\int_0^{1/2} 
\xi^{2k}e^{-2\xi^2}
\,d\xi
\right)^{\frac{1}{2}}
\left|
\sum_{j=0}^{k/2}
(-1)^j M_{2j}(v)
\right|
}, 
& k\equiv0\pmod{2}, \\[22pt]
\displaystyle{
\left(
2\int_0^{1/2} 
\xi^{2k}e^{-2\xi^2}
\,d\xi
\right)^{\frac{1}{2}}
\left|
\sum_{j=0}^{(k-1)/2}
(-1)^j M_{2j+1}(v)
\right|
}, 
& k\equiv1\pmod{2}.
\end{cases}
\end{align*}
\end{proposition}

\begin{proposition}
\label{prop:2}
Let $n\ge2$ and $k\in\bf{N}_0$. 
For $v\in L^{1,k}({\bf R}^n)$, it holds that 
\begin{align*}
\biggr\|B_k^v e^{-|\xi|^2}\biggr\|_{L^2(|\xi|\le1/2)}=
\begin{cases}
\displaystyle{
\left(
\int_{|\xi|\le1/2} 
e^{-2|\xi|^2}
\,d\xi
\right)^{\frac{1}{2}}
|M_0(v)|
}, 
& k=0, \\[22pt]
\displaystyle{
\left(
\int_{|\xi|\le1/2} 
\xi_1^2 e^{-2|\xi|^2}
\,d\xi
\right)^{\frac{1}{2}}
\left(
\sum_{|\alpha|=1}
M_\alpha(v)^2
\right)^\frac{1}{2}
}, 
& k=1.
\end{cases}
\end{align*}
Furthermore, 
\begin{align*}
\biggr\|B_2^v e^{-|\xi|^2}\biggr\|_{L^2(|\xi|\le1/2)}
=\left[
C_1\sum_{j=1}^n 
V_j^2 
+C_{12}\sum_{1\le j<k\le n} 
\biggr(2V_j V_k +W_{j,k}^2\biggr)
\right]^\frac{1}{2}
\end{align*}
where 
\[
C_1:=
\int_{|\xi|\le1/2}
\xi_1^4 e^{-2|\xi|^2}\,d\xi,
\qquad 
C_{12}:=
\int_{|\xi|\le1/2}
\xi_1^2 \xi_2^2 e^{-2|\xi|^2}\,d\xi,
\]
\[
V_j
:=\int_{{\bf R}^{n}}v(x)\,dx
-\frac{1}{2}\int_{{\bf R}^{n}}x_j^2 v(x)\,dx,
\qquad 
W_{k,j}
:=\int_{{\bf R}^{n}}x_j x_k v(x)\,dx.
\]
\end{proposition}

\section{Proofs}

In this section, let us prove our main results. We first present the following generalized lemma. This lemma with $\gamma = 1$ has its origin in \cite[Lemma 3.1]{I-0}.
\begin{lemma}
\label{lem:5.1}
Let $n\ge1$ and $\gamma\ge0$. 
Then it holds  
\begin{align}
\label{5.1}
\left|
\hat{f}(\xi)
-\sum_{|\alpha|\le[\gamma]}
M_\alpha(f)
(i\xi)^\alpha
\right|
\le C_\gamma |\xi|^\gamma
\int_{{\bf R}^n} |x|^\gamma |f(x)|\,dx
\end{align}
for $f\in L^{1,\gamma}({\bf R}^n)$. 
Here $C_\gamma>0$ is a constant independent of $\xi$ and $f$. 
\end{lemma}
{\bf Proof.} 
By the Taylor theorem we see that 
\begin{align*}
e^{-ix\cdot\xi}
=\sum_{|\alpha|\le[\gamma]}
\frac{(-1)^{|\alpha|}}{\alpha!}x^\alpha (i\xi)^\alpha
+\mathcal{R}_{[\gamma]+1}(x,\xi),
\end{align*}
where 
\begin{align*}
\mathcal{R}_{[\gamma]+1}(x,\xi)
& :=\frac{1}{[\gamma]!}
\int_0^1 (1-\tau)^{[\gamma]}
\frac{d^{[\gamma]+1}} {d\tau^{[\gamma]+1}} e^{-i\tau x\cdot\xi}\,
d\tau \\
& =\frac{(-i)^{[\gamma]+1}}{[\gamma]!}
(x\cdot\xi)^{[\gamma]+1}
\int_0^1 (1-\tau)^{[\gamma]}
e^{-i\tau x\cdot\xi}\,
d\tau.
\end{align*}
Now, it follows that 
\begin{align*}
C_\gamma
:=\sup_{x,\,\xi\in{\bf R}^n\setminus\{0\}}
\frac{|\mathcal{R}_{[\gamma]+1}(x,\xi)|}{|x|^\gamma |\xi|^\gamma}
<+\infty,
\end{align*}
since 
\[
\sup_{0<|x| |\xi|\le1}
\frac{|\mathcal{R}_{[\gamma]+1}(x,\xi)|}{|x|^\gamma |\xi|^\gamma}
\le
\frac{1}{[\gamma]!}
\sup_{0<|x| |\xi|\le1}
(|x| |\xi|)^{[\gamma]+1-\gamma}
\le \frac{1}{[\gamma]!}, 
\]
\[
\sup_{|x| |\xi|\ge1}
\frac{|\mathcal{R}_{[\gamma]+1}(x,\xi)|}{|x|^\gamma |\xi|^\gamma}
\le
\sup_{|x| |\xi|\ge1}
\left(
\frac{1}{(|x| |\xi|)^\gamma}
+\sum_{|\alpha|\le[\gamma]}
\frac{1}{\alpha!}
\frac{1}{(|x| |\xi|)^{\gamma-|\alpha|}}
\right)
\le 1+\sum_{|\alpha|\le[\gamma]}
\frac{1}{\alpha!}.
\]
Therefore we obtain  
\begin{align*}
\left|
\hat{f}(\xi)
-\sum_{|\alpha|\le[\gamma]}
\frac{(-1)^{|\alpha|}}{\alpha!}
\left(
\int_{{\bf R}^n} 
x^\alpha
f(x)\,
dx
\right)
(i\xi)^\alpha
\right| 
& \le \int_{{\bf R}^n} 
\big|\mathcal{R}_{[\gamma]+1}(x,\xi) v(x)\big| \, dx \\
& \le C_\gamma |\xi|^\gamma
\int_{{\bf R}^n} |x|^\gamma |v(x)|\,dx,
\end{align*}
which implies \eqref{5.1}. 
$\Box$

Furthermore, we prepare the following key lemma, which will be essentially used to obtain a higher order asymptotic expansion of the solution to \eqref{1.1}. This can be derived with the help of Lemma~\ref{lem:5.1} just proved. 
\begin{lemma}
\label{lem:5.2}
Let $n\ge1$ and $\gamma\ge0$. If $v\in L^2({\bf R}^n)\cap L^{1,\gamma}({\bf R}^n)$, 
then it holds 
\begin{align}
\label{5.2}
\big|F^v(\xi)-A_{[\gamma]}^v(\xi)\big|
\le C|\xi|^\gamma \|v\|_{1,\gamma}
\end{align}
for $\xi\in{\bf R}^n$ with $|\xi|\le1/2$. 
Here $C>0$ is a constant independent of $\xi$ and $v$. 
\end{lemma}
{\bf Proof.} 
If $[\gamma]$ is even, 
it holds that 
\begin{align*}
\frac{1}{1-|\xi|^2}
=\sum_{j=0}^{[\gamma]/2}
|\xi|^{2j}
+\frac{|\xi|^{[\gamma]+2}}{1-|\xi|^2}, 
\end{align*}
\begin{align*}
\hat{v}(\xi)
=\sum_{|\alpha|\le[\gamma]}
M_\alpha(v)(i\xi)^\alpha 
+\int_{{\bf R}^n}\mathcal{R}_{[\gamma]+1}(x,\xi)v(x)\,dx. 
\end{align*}
Thus it follows that 
\begin{align*}
F^v(\xi)
& =\left(
\sum_{j=0}^{[\gamma]/2} |\xi|^{2j} 
+\frac{|\xi|^{[\gamma]+2}}{1-|\xi|^2}
\right)
\left(
\sum_{|\alpha|\le[\gamma]}
M_\alpha(v)(i\xi)^\alpha 
+\int_{{\bf R}^n}\mathcal{R}_{[\gamma]+1}(x,\xi)v(x)\,dx
\right) \\
& =\left(
\sum_{j=0}^{[\gamma]/2} |\xi|^{2j} 
\right)
\left(
\sum_{|\alpha|\le[\gamma]}
M_\alpha(v)(i\xi)^\alpha 
\right) \\
& \qquad\qquad
+\left(
\sum_{j=0}^{[\gamma]/2} |\xi|^{2j} 
\right)
\int_{{\bf R}^n}\mathcal{R}_{[\gamma]+1}(x,\xi)v(x)\,dx 
+\frac{|\xi|^{[\gamma]+2}}{1-|\xi|^2} \hat{v}(\xi). 
\end{align*}
Furthermore, we see that 
\begin{align*}
& \left(
\sum_{j=0}^{[\gamma]/2} |\xi|^{2j} 
\right)
\left(
\sum_{|\alpha|\le[\gamma]}
M_\alpha(v)(i\xi)^\alpha 
\right)
=\left(
\sum_{j=0}^{[\gamma]/2} |\xi|^{[\gamma]-2j} 
\right)
\left(
\sum_{|\alpha|\le[\gamma]}
M_\alpha(v)(i\xi)^\alpha 
\right) \\
& =\left(
\sum_{j=0}^{[\gamma]/2} |\xi|^{[\gamma]-2j} 
\right)
\left(
\sum_{|\alpha|\le2j}
M_\alpha(v)(i\xi)^\alpha
+\sum_{2j+1\le|\alpha|\le[\gamma]}
M_\alpha(v)(i\xi)^\alpha 
\right) \\
& =A_{[\gamma]}^v(\xi)
+\sum_{j=0}^{[\gamma]/2-1}
|\xi|^{[\gamma]-2j}
\left(
\sum_{2j+1\le|\alpha|\le[\gamma]}
M_\alpha(v)(i\xi)^\alpha 
\right).
\end{align*}
Note that when $[\gamma]=0$, 
we regard the sum in the second term of the right-hand side above as zero. 
So we arrive at 
\begin{align*}
F^v(\xi)-A_{[\gamma]}^v(\xi)
& =\sum_{j=0}^{[\gamma]/2-1}
|\xi|^{[\gamma]-2j}
\left(
\sum_{2j+1\le|\alpha|\le[\gamma]}
M_\alpha(v)(i\xi)^\alpha 
\right) \\
& \qquad\qquad
+\left(
\sum_{j=0}^{[\gamma]/2} |\xi|^{2j} 
\right)
\int_{{\bf R}^n}\mathcal{R}_{[\gamma]+1}(x,\xi)v(x)\,dx 
+\frac{|\xi|^{[\gamma]+2}}{1-|\xi|^2} \hat{v}(\xi),
\end{align*}
which implies 
\begin{align*}
\big|F^v(\xi)-A_{[\gamma]}^v(\xi)\big|
& \le C|\xi|^{[\gamma]+1} \|v\|_{1,[\gamma]} 
+\int_{{\bf R}^n}\big|\mathcal{R}_{[\gamma]+1}(x,\xi)v(x)\big|\,dx 
+C|\xi|^{[\gamma]+2} |\hat{v}(\xi)| \\
& \le C|\xi|^{[\gamma]+1} \|v\|_{1,[\gamma]} 
+C |\xi|^\gamma \|v\|_{1,\gamma}
+C|\xi|^{[\gamma]+2} \|v\|_1 \\
& \le C |\xi|^\gamma \|v\|_{1,\gamma}
\end{align*}
for $\xi\in{\bf R}^n$ with $|\xi|\le1/2$. 
At the second inequality we use \eqref{5.1} in Lemma~\ref{lem:5.1}. 
When $[\gamma]$ is odd, 
the proof is similar. 
$\Box$

Now, let us prove Theorem~\ref{thm:1} by using Lemma~\ref{lem:5.2}.\\

{\bf Proof of Theorem~\ref{thm:1}.} 
For simplicity, we write $v:=u_0+u_1$. 
We first establish the desired estimate in the low frequency region for $\xi$ with $|\xi|\le1/2$. Indeed, from \eqref{2.1} and \eqref{5.2} in Lemma~\ref{lem:5.2}, 
it follows that 
\begin{align*}
\biggr\|
\hat{u}(t)
-A_{[\gamma]}^v e^{-t|\xi|^2}
\biggr\|_{L^2(|\xi|\le1/2)}
& \le C\biggr\|
|\xi|^\gamma e^{-t|\xi|^2}
\biggr\|_{L^2(|\xi|\le1/2)}
\|v\|_{1,\gamma}
+Ce^{-t}\big(\|u_0\|_2+\|u_1\|_2\big) \\
& \le C(1+t)^{-\frac{n}{4}-\frac{\gamma}{2}}
\|v\|_{1,\gamma}
+Ce^{-t}\big(\|u_0\|_2+\|u_1\|_2\big)
\end{align*}
for $t>0$. 

Next, we consider the middle frequency part for $\xi$ satisfying $1/2\le|\xi|\le2$. For this purpose we use the expression \eqref{2.2}. 
By the Taylor theorem we have 
\begin{align*}
\frac{e^{-t|\xi|^2}-e^{-t}}{1-|\xi|^2}
=-\frac{e^{-t|\xi|^2}-e^{-t}}{|\xi|^2-1} 
=-\left.\frac{\partial}{\partial r}e^{-tr}\right|_{r=\tau}
=te^{-t\tau}
\end{align*}
for some $1/2\le\tau\le2$. 
So we have 
\[
\sup_{1/2\le|\xi|\le2}
\left|
\frac{e^{-t|\xi|^2}-e^{-t}}{1-|\xi|^2}
\right|
\le Ce^{-ct}
\]
for $t>0$. 
Furthermore, it follows from 
\begin{align*}
\frac{e^{-t|\xi|^2}-|\xi|^2 e^{-t}}{1-|\xi|^2}
=-\frac{e^{-t|\xi|^2}-e^{-t}}{|\xi|^2-1}+e^{-t}
\end{align*}
that 
\[
\sup_{1/2\le|\xi|\le2}
\left|
\frac{e^{-t|\xi|^2}-|\xi|^2 e^{-t}}{1-|\xi|^2}
\right|
\le Ce^{-ct}
\]
for $t>0$. 
Thus it holds 
\begin{align*}
\|\hat{u}(t)\|_{L^2(1/2\le|\xi|\le2)}
\le Ce^{-ct}\big(\|u_0\|_2+\|u_1\|_2\big)
\end{align*}
for $t>0$. 
Here $c>0$ is a constant independent of $t$, $u_0$ and $u_1$. 

Finally, let us obtain the desired estimate in the high frequency region for $\xi$ satisfying $|\xi|\ge2$. From \eqref{2.3} one can easily find that 
\begin{align*}
\|\hat{u}(t)\|_{L^2(|\xi|\ge2)}
\le Ce^{-t}\big(\|u_0\|_2+\|u_1\|_2\big)
\end{align*}
for $t>0$. 
Here we have just used the fact:
\[
\sup_{|\xi|\ge2} \frac{|\xi|^2}{|\xi|^2-1}<\infty.
\]
On the other hand, we have 
\begin{align*}
\biggr\|A_{[\gamma]}^v e^{-t|\xi|^2}\biggr\|_{L^2(|\xi|\ge1/2)}
& \le C\biggr\||\xi|^{[\gamma]}e^{-t|\xi|^2}\biggr\|_{L^2(|\xi|\ge1/2)}
\|v\|_{1,[\gamma]} \\
& \le Ce^{-ct}\|v\|_{1,[\gamma]}
\end{align*}
for $t\ge1$. The desired estimate now follows. 
$\Box$

Next, we give the following lemma on the solution for the heat equation (see also the related results in Appendix). 
Lemma~\ref{lem:5.3} leads to Lemma~\ref{lem:5.4} which is a key to the proof of Theorem~\ref{thm:2}. 
\begin{lemma}
\label{lem:5.3}
Let $n\ge1$, $\ell\ge0$ and 
$v\in L^{1,\gamma}({\bf R}^n)$ with $\gamma\ge0$. 
Then it holds that 
\begin{align}
\label{5.3}
\lim_{t\to\infty}
t^{\frac{n}{4}+\frac{\gamma}{2}+\frac{\ell}{2}}
\left\|
|\xi|^{\ell}
\left(
e^{-t|\xi|^2}\hat{v}
-\sum_{|\alpha|\le [\gamma]}
M_\alpha(v) (i\xi)^\alpha
e^{-t|\xi|^2}
\right)
\right\|_2
=0.
\end{align}
\end{lemma}
{\bf Proof.} 
We have already seen that 
\begin{align*}
\hat{v}(\xi)
-\sum_{|\alpha|\le[\gamma]}
M_\alpha(v) (i\xi)^\alpha 
& =\int_{{\bf R}^n}
\mathcal{R}_{[\gamma]+1}(x,\xi)v(x)\,dx
\end{align*}
in the proof of Lemma~\ref{lem:5.1}. 
By the form of $\mathcal{R}_{[\gamma]+1}$, 
there exists a constant $C>0$ such that 
\[
|\mathcal{R}_{[\gamma]+1}(x,\xi)|
\le C(|x||\xi|)^\gamma
\]
for $x\in{\bf R}^n$ and $\xi\in{\bf R}^n$. 
Now we calculate 
\begin{align*}
& \left\|
|\xi|^{\ell}
\left(
e^{-t|\xi|^2}\hat{v}
-\sum_{|\alpha|\le[\gamma]}
M_\alpha(v) (i\xi)^\alpha
e^{-t|\xi|^2}
\right)
\right\|_2^2 \\
& =
\int_{{\bf R}^n}
|\xi|^{2\ell+2\gamma}
\left|
\int_{{\bf R}^n}
\frac{\mathcal{R}_{[\gamma]+1}(x,\xi)}{|x|^\gamma |\xi|^\gamma}
|x|^\gamma v(x)\,dx
\right|^2
e^{-2t|\xi|^2}
\,d\xi \\
& =t^{-\frac{n}{2}-\ell-\gamma}
\int_{{\bf R}^n}
|\eta|^{2\ell+2\gamma}
\left|
\int_{{\bf R}^n}
\frac{\mathcal{R}_{[\gamma]+1}(x,\eta/\sqrt{t})}{|x|^\gamma |\eta/\sqrt{t}|^\gamma}
|x|^\gamma v(x)\,dx
\right|^2
e^{-2|\eta|^2}
\,d\eta.
\end{align*}
In the proof of Lemma~\ref{lem:5.1}, 
we have also proved that 
\[
\sup_{x,\,\xi\in{\bf R}^n\setminus\{0\}}
\frac{|\mathcal{R}_{[\gamma]+1}(x,\xi)|}{|x|^\gamma |\xi|^\gamma}
<\infty. 
\]
Furthermore, we can easily see that 
\[
\lim_{\xi\to0}
\frac{\mathcal{R}_{[\gamma]+1}(x,\xi)}{|x|^\gamma |\xi|^\gamma}
=0
\] 
for fixed $x\in{\bf R}^n\setminus\{0\}$. 
Thus Lebesgue's dominated convergence theorem yields 
\begin{align*}
& \lim_{t\to\infty}
t^{\frac{n}{2}+\ell+\gamma}
\left\|
|\xi|^{\ell}
\left(
e^{-t|\xi|^2}\hat{v}
-\sum_{|\alpha|\le[\gamma]}
M_\alpha(v) (i\xi)^\alpha
e^{-t|\xi|^2}
\right)
\right\|_2^2 \\
& =
\int_{{\bf R}^n}
|\eta|^{2\ell+2\gamma}
\left|
\int_{{\bf R}^n}
\lim_{t\to\infty}
\frac{\mathcal{R}_{[\gamma]+1}(x,\eta/\sqrt{t})}{|x|^\gamma |\eta/\sqrt{t}|^\gamma}
|x|^\gamma v(x)\,dx
\right|^2
e^{-2|\eta|^2}
\,d\eta 
=0,
\end{align*}
which implies the desired equality. 
$\Box$

\begin{lemma}
\label{lem:5.4}
Let $n\ge1$, $\ell\ge0$ and $v\in L^{1,k}({\bf R}^n)$ with $k\in\bf{N}_0$. 
Then it holds that 
\begin{align}
\label{5.4}
\lim_{t\to\infty}
t^{\frac{n}{4}+\frac{k}{2}+\frac{\ell}{2}}
\left\|
|\xi|^\ell
\left(
F^v e^{-t|\xi|^2}
-A_k^v e^{-t|\xi|^2}
\right)
\right\|_{L^2(|\xi|\le1/2)}
=0.
\end{align}
\end{lemma}
{\bf Proof.} 
Take any $k\in{\bf N}_0$ with $k\equiv0\pmod{2}$. 
By the definitions of $F^v$ and $A_k^v$, 
we have 
\begin{align}
\notag
& F^v(\xi)-A_k^v(\xi)
=\frac{1}{1-|\xi|^2}\hat{v}(\xi)
-\sum_{j=0}^{k/2}
\left(
|\xi|^{k-2j}
\sum_{|\alpha|\le2j}
M_\alpha(v)(i\xi)^\alpha
\right) \\
\notag
& =\left(
\sum_{j=0}^{k/2}|\xi|^{2j}
+\frac{|\xi|^{k+2}}{1-|\xi|^2}\right)\hat{v}(\xi)
-\sum_{j=0}^{k/2}
\left(
|\xi|^{k-2j}
\sum_{|\alpha|\le2j}
M_\alpha(v)(i\xi)^\alpha
\right) \\
\notag
& =\left(
\sum_{j=0}^{k/2}|\xi|^{k-2j}
+\frac{|\xi|^{k+2}}{1-|\xi|^2}\right)\hat{v}(\xi)
-\sum_{j=0}^{k/2}
\left(
|\xi|^{k-2j}
\sum_{|\alpha|\le2j}
M_\alpha(v)(i\xi)^\alpha
\right) \\
\label{5.5}
& =\sum_{j=0}^{k/2}
|\xi|^{k-2j}
\left(
\hat{v}(\xi)
-\sum_{|\alpha|\le2j}
M_\alpha(v)(i\xi)^\alpha
\right)
+\frac{|\xi|^{k+2}}{1-|\xi|^2}\hat{v}(\xi).
\end{align}
With the aid of \eqref{5.3} in Lemma~\ref{lem:5.3}, 
one can derive 
\begin{align*}
& \left\|
|\xi|^{\ell+k-2j}\left(
e^{-t|\xi|^2}\hat{v}(\xi)
-\sum_{|\alpha|\le 2j}M_\alpha(v)(i\xi)^\alpha e^{-t|\xi|^2}
\right)
\right\|_{L^{2}(|\xi|\le1/2)} \\
& \le \left\|
|\xi|^{\ell+k-2j}\left(
e^{-t|\xi|^2}\hat{v}(\xi)
-\sum_{|\alpha|\le 2j}M_\alpha(v)(i\xi)^\alpha e^{-t|\xi|^2}
\right)
\right\|_2
=o(t^{-\frac{n}{4}-\frac{k}{2}-\frac{\ell}{2}})
\end{align*}
for each $j=0,1,\cdots,k/2$. 
Moreover, it is easy to check that
\[
\left\| 
|\xi|^{\ell}\frac{|\xi|^{k+2}}{1-|\xi|^{2}}
e^{-t|\xi|^{2}}\hat{v}
\right\|_{L^{2}(|\xi|\le1/2)} 
=O(t^{-\frac{n}{4}-\frac{k}{2}-\frac{\ell}{2}-1}), 
\qquad
t \to \infty.
\]
This shows the validity of the statement for $k \equiv 0$ (mod $2$). The other case is also similar. 
$\Box$
\\

{\bf Proof of Theorem~\ref{thm:2}.} 
Recalling property~(A) and the solution formula~\eqref{2.1}, we have 
\begin{align*}
& \biggr\|\hat{u}(t)-A_{k-1}^v e^{-t|\xi|^2}\biggr\|_2 \\
& \ge\biggr\|
F^v e^{-t|\xi|^2}
-A_{k-1}^v e^{-t|\xi|^2}
\biggr\|_{L^2(|\xi|\le1/2)} 
-o(e^{-ct}) \\
& =\biggr\|
\left(
F^v e^{-t|\xi|^2}
-A_k^v e^{-t|\xi|^2}
\right)
+B_k^v e^{-t|\xi|^2}
\biggr\|_{L^2(|\xi|\le1/2)} 
-o(e^{-ct}) \\
& \ge \left\|
B_k^v e^{-t|\xi|^2}
\right\|_{L^2(|\xi|\le1/2)} 
-\biggr\|
F^v e^{-t|\xi|^2}
-A_k^v e^{-t|\xi|^2}
\biggr\|_{L^2(|\xi|\le1/2)} 
-o(e^{-ct})
\end{align*}
as $t\to\infty$. 
By \eqref{5.4} in Lemma~\ref{lem:5.4}, one has 
\[
\biggr\|
F^v e^{-t|\xi|^2}
-A_k^v e^{-t|\xi|^2}
\biggr\|_{L^2(|\xi|\le1/2)}
=o(t^{-\frac{n}{4}-\frac{k}{2}}), 
\qquad
t\to\infty. 
\]
Furthermore, by changing variables such that $\eta=\sqrt{t}\xi$ we see 
\begin{align*}
\left\|
B_k^v(\xi) e^{-t|\xi|^2}
\right\|_{L^2(|\xi|\le1/2)}
& =t^{-\frac{n}{4}-\frac{k}{2}}
\left\|
B_k^v(\eta) e^{-|\eta|^2}
\right\|_{L^2(|\eta|\le\sqrt{t}/2)} \\
& \ge t^{-\frac{n}{4}-\frac{k}{2}}
\left\|
B_k^v(\eta) e^{-|\eta|^2}
\right\|_{L^2(|\eta|\le1/2)}, 
\qquad 
t\ge1. 
\end{align*}
Here we have just used the property~(C). Therefore the proof is now complete. 
$\Box$
\\

{\bf Proof of Proposition~\ref{prop:1}.} 
By definition of $B_k^v$, we have 
\begin{align*}
B_k^v(\xi):=
\begin{cases}
\displaystyle{
\left(
\sum_{j=0}^{k/2}
(-1)^j M_{2j}(v)
\right)
\xi^k
},
& k\equiv0\pmod{2},\\[22pt]
\displaystyle{
\left(
i \sum_{j=0}^{(k-1)/2}
(-1)^j M_{2j+1}(v)
\right)
\xi^k
},
& k\equiv1\pmod{2}, 
\end{cases}
\end{align*}
and thus the statement of proposition now follows. 
$\Box$

\begin{remark}
{\rm Further calculation shows 
\[
\sum_{j=0}^{k/2}
(-1)^j M_{2j}(v)
=\int_{-\infty}^\infty 
\left(
\sum_{j=0}^{k/2}
\frac{(-1)^j}{(2j)!}
x^{2j}
\right)
v(x)\,dx,
\qquad
k\equiv0\pmod{2}, 
\]
\[
i\sum_{j=0}^{(k-1)/2}
(-1)^j M_{2j+1}(v)
=-i \int_{-\infty}^\infty 
\left(
\sum_{j=0}^{(k-1)/2}
\frac{(-1)^j}{(2j+1)!}
x^{2j+1}
\right)
v(x)\,dx,
\qquad
k\equiv1\pmod{2}. 
\]}
\end{remark}
\vspace{0.2cm}

{\bf Proof of Proposition~\ref{prop:2}.} 
First, we consider the case $k=0$. 
However, this case is trivial since $B_0^v(\xi)\equiv M_0(v)$. 
Next, we deal with the case $k = 1,2$. 
Let $v\in L^{1,1}({\bf R}^n)$. 
For simplicity, 
using the symbols of the standard basis $\{e_j\}_{j=1,\dots,n}$ for ${\bf R}^n$, 
we write 
\[
M_j(v):=M_{e_j}(v), 
\qquad
j=1,2,\dots,n. 
\]
Then we see that 
\begin{align*}
\big|B_1^v(\xi)\big|^2
& =\left|\sum_{|\alpha|=1}M_\alpha(v)(i\xi)^\alpha \right|^2 
=\left(\sum_{|\alpha|=1}M_\alpha(v)\xi^\alpha \right)^2 \\
& =\sum_{j=1}^n 
M_j(v)^2 \xi_j^2
+2\sum_{1\le j<k\le n} 
M_j(v) M_k(v) \xi_j \xi_k. 
\end{align*}
Note that 
\[
\int_{|\xi|\le1/2}
\xi_j \xi_k e^{-2|\xi|^2}
\,d\xi
=0, 
\qquad
1\le j<k\le n,
\]
and so we have 
\begin{align*}
\left\|
B_1^v e^{-|\xi|^2}
\right\|_{L^2(|\xi|\le1/2)}^2 
=\left(
\int_{|\xi|\le1/2}
\xi_1^2 e^{-2|\xi|^2}
\,d\xi
\right)
\sum_{j=1}^n M_j(v)^2. 
\end{align*}
Thus we obtain the proposition in the case $k=1$. 

Finally, for $v\in L^{1,2}({\bf R}^n)$, 
we rearrange  
\begin{align*}
B_2^v(\xi)
& =|\xi|^2 M_0(v)
+\sum_{|\alpha|=2}
M_\alpha(v)(i\xi)^\alpha \\
& =\sum_{j=1}^n 
\left(\int_{{\bf R}^n}v(x)\,dx\right)
\xi_j^2 
-\frac{1}{2} 
\sum_{j=1}^n 
\left(\int_{{\bf R}^n}x_j^2 v(x)\,dx\right)
\xi_j^2 
-\sum_{1\le j<k\le n} 
\left(\int_{{\bf R}^n}x_j x_k v(x)\,dx\right)
\xi_j \xi_k \\
& =\sum_{j=1}^n 
V_j
\xi_j^2 
-\sum_{1\le j<k\le n} 
W_{j,k}
\xi_j \xi_k.
\end{align*}
So we arrive at 
\begin{align*}
B_2^v(\xi)^2
=\sum_{j=1}^n V_j^2 \xi_j^4
+\sum_{1\le j<k\le n} 
\biggr(2V_j V_k +W_{j,k}^2\biggr)
\xi_j^2 \xi_k^2 
+I(\xi). 
\end{align*}
The function $I(\xi)$ consists of odd functions for some $\xi_j$ and so 
\[
\int_{|\xi|\le1/2}
I(\xi)e^{-2|\xi|^2}\,d\xi
=0.
\]
Thus we obtain  
\begin{align*}
\notag
& \biggr\|B_2^v e^{-|\xi|^2}\biggr\|_{L^2(|\xi|\le1/2)}^2 \\
\notag
& =\sum_{j=1}^n 
V_j^2 
\int_{|\xi|\le1/2}
\xi_j^4 e^{-2|\xi|^2}\,d\xi
+\sum_{1\le j<k\le n} 
\biggr(2V_j V_k +W_{j,k}^2\biggr)
\int_{|\xi|\le1/2}
\xi_j^2 \xi_k^2 e^{-2|\xi|^2}\,d\xi \\
\label{B2v}
& =C_1\sum_{j=1}^n 
V_j^2 
+C_{12} \sum_{1\le j<k\le n} 
\biggr(2V_j V_k +W_{j,k}^2\biggr),
\end{align*}
and the proof is now complete. 
$\Box$ 
\vspace{0.2cm}
\begin{example}
\label{ex:5.1}
{\rm Here we see that the positivity of $\biggr\|B_k^v e^{-|\xi|^2}\biggr\|_{L^2(|\xi|\le1/2)}$, i.e., the non-zero property of $B_k^v e^{-|\xi|^2}$ is a delicate problem. 
On this interest see also Proposition~\ref{prop:6.2} below. 

Let $v\in L^{1,k}({\bf R}^n)$ with $k\in{\bf N}_0$. 
If $M_\alpha(v)=0$ for all $|\alpha|\le k$, 
it is obvious that $B_k^v(\xi)\equiv0$. 
The contrary is also true for $k=0,1$, but this shall not be true when $k\ge2$. 
Actually we can find $v\in L^{1,2}({\bf R}^n)$ satisfying $M_\alpha(v)\not=0$ for some $\alpha\in{\bf N}_0^n$ with $|\alpha|\le2$ and  $B_2^v(\xi)\equiv0$. 
A typical example is 
\[
v(x):=\exp\left(-\frac{|x|^2}{4}\right).
\] 
Then $v\in L^{1,2}({\bf R}^n)$. 
If $n=1$, we can easily obtain 
\[
M_0(v)=\int_{-\infty}^\infty v(x)\,dx
=2\sqrt{\pi}, 
\qquad
M_2(v)=\int_{-\infty}^\infty x^2 v(x)\,dx
=4\sqrt{\pi},
\]
and thus $B_2^v(\xi)\equiv0$. 
From the above calculation, in the case $n\ge2$, $v$ satisfies 
\begin{itemize}
\item[\rm(i)] $M_0(v)=(2\sqrt{\pi})^n$ 
and 
$M_{2e_j}(v)=2(2\sqrt{\pi})^n$ for all $j=1,\cdots,n${\rm ;}
\item[\rm(ii)]
$V_j=0$ for all $j=1,\dots,n${\rm ;}   
\item[\rm(iii)]
$W_{k,j}=0$ for all $1\le j<k\le n$. 
\end{itemize}
Assertions {\rm (ii)} and {\rm (iii)} imply $B_2^v(\xi)\equiv0$.}
\end{example}

\section{Heat equation}
In this section, we re-study the heat equation: 
\begin{align}
\label{6.1}
\begin{cases}
u_t-\Delta u=0, & t>0,\quad x\in{\bf R}^n, \\
u(0,x)=v(x), & x\in{\bf R}^n,
\end{cases}
\end{align}
where $n\ge1$ and $v\in L^{1,k}({\bf R}^n)$ with $k\in{\bf N}_0$. 
Now we state some results corresponding to Theorem~\ref{thm:2}, Propositions~\ref{prop:1} and \ref{prop:2} in order to compare them. 

\begin{Proposition}
\label{prop:6.1}
Let $n\ge1$ and $k\in\bf{N}_0$, and let $\hat{u}$ be the solution to \eqref{6.1} with $v\in L^2({\bf R}^{n})\cap L^{1,k}({\bf R}^{n})$, i.e., $\hat{u}(t,\xi)=e^{-t|\xi|^2}\hat{v}$. 
Then there exists a constant $\delta\ge1$ depending on $v$ such that
\begin{align*}
\biggr\|
\hat{u}(t)
-\sum_{|\alpha|\le k-1}
M_\alpha(v)(i\xi)^\alpha 
e^{-t|\xi|^2}
\biggr\|_2
\ge \frac{1}{2}
\biggr\|
\sum_{|\alpha|=k}
M_\alpha(v)(i\xi)^\alpha 
e^{-|\xi|^2}
\biggr\|_2
t^{-\frac{n}{4}-\frac{k}{2}}
\end{align*}
for $t\ge\delta$.
\end{Proposition}

Proposition~\ref{prop:6.1} implies that 
\[
C_k^v(\xi)
:=\sum_{|\alpha|=k}
M_\alpha(v)(i\xi)^\alpha 
\quad
\mbox{corresponds to}
\quad
B_k^v(\xi)
=|\xi|^2 B_{k-2}^v(\xi)
+\sum_{|\alpha|=k}
M_\alpha(v)(i\xi)^\alpha
\]
in the heat flow case 
(see \eqref{4.2} with property~(B)). 
In this sense $B_k^v$ appearing in \eqref{4.2} describes a difference between the solution to \eqref{1.1} and that of \eqref{6.1} although $B_k^v(\xi)\equiv C_k^v(\xi)$ for $k=0,1$. 
The case of $k\ge2$ is meaningful in this paper.
In addition, the following proposition shows that 
the function $C_k^v$ possesses a simple structure than $B_k^v$. 
\begin{proposition}
\label{prop:6.2}
Let $n\ge1$ and $v\in L^{1,k}({\bf R}^n)$ with $k\in\bf{N}_0$. 
Then the followings are equivalent{\rm :}
\begin{itemize}
\item[\rm{(1)}]
$C_k^v(\xi)\equiv0${\rm ;}
\item[\rm{(2)}]
$M_\alpha(v)=0$ for all $|\alpha|=k$.
\end{itemize}
Furthermore, it holds that 
\begin{align*}
\biggr\|
C_k^v e^{-|\xi|^2}
\biggr\|_2^2
=\sum_{|\alpha|=k}
\left[
\left(
\int_{\bf{R}^n}
\xi^{2\alpha} e^{-2|\xi|^2}
\, d\xi
\right)
\sum_{\substack{\beta_1+\beta_2=2\alpha, \\[3pt] 
|\beta_1|=|\beta_2|=k}}
M_{\beta_1}(v)M_{\beta_2}(v)
\right].
\end{align*}
\end{proposition}
\section{Appendix}
In this section we shall review the results from \cite{IKM}, which dealt with the fractional heat equation. 
\begin{theorem}{\rm (\cite{IKM})}
Let $\gamma\ge0$, $\beta\in{\bf N}_{0}^{n}$ and $v\in L^{1,\gamma}({\bf R}^n)$. 
Then 
\begin{equation}
\label{7.1}
\begin{split}
& \left\|
\partial_x^\beta
\left[
\int_{{\bf R}^n}
G(t,x-y)v(y)\,dy
-\sum_{|\alpha|\le[\gamma]}
M_\alpha(v)
(\partial_x^\alpha G)(t,x)
\right]
\right\|_q \\
& \hspace{5cm}
\le Ct^{-\frac{n}{2}(1-\frac{1}{q})-\frac{\gamma}{2}-\frac{|\beta|}{2}}
\|v\|_{1,\gamma}, 
\qquad 
t>0, 
\end{split}
\end{equation}
for $1\le q\le\infty$. 
Here $C>0$ is a constant independent of $t$ and $v$. 
Furthermore, for any $1\le q\le\infty$ it holds that 
\begin{align}
\label{7.2}
\lim_{t\to\infty}
t^{\frac{n}{2}(1-\frac{1}{q})+\frac{\gamma}{2}+\frac{|\beta|}{2}}
\left\|
\partial_x^\beta
\left[
\int_{{\bf R}^n}
G(t,x-y)v(y)\,dy
-\sum_{|\alpha|\le[\gamma]}
M_\alpha(v)
(\partial_x^\alpha G)(t,x)
\right]
\right\|_q
=0. 
\end{align}
\end{theorem}
\begin{Remark}
{\rm In the case when $\gamma = 0$ and $\beta = 0$ in (7.2), the inequality has already been introduced in \cite[Lemma 3.2]{K} (see also \cite{DZ} for a more general result).}
\end{Remark}

As a result of Theorem 7.1, one can get the following statements.\\ 
If $\ell\ge0$ and $v\in L^{1,\gamma}({\bf R}^n)$ with $\gamma\ge0$, 
it follows from \eqref{7.1} with $q=2$ that 
\begin{align}
\label{7.3}
\left\|
|\xi|^\ell
\left(
e^{-t|\xi|^2}
\hat{v}
-\sum_{|\alpha|\le[\gamma]}
M_\alpha(v)
(i\xi)^\alpha e^{-t|\xi|^2}
\right)
\right\|_2
\le Ct^{-\frac{n}{4}-\frac{\gamma}{2}-\frac{\ell}{2}}
\|v\|_{1,\gamma}, 
\qquad 
t>0, 
\end{align}
and \eqref{7.2} with $q=2$ is equivalent to
\begin{align}
\tag{5.3}
\lim_{t\to\infty}
t^{\frac{n}{4}+\frac{\gamma}{2}+\frac{\ell}{2}}
\left\|
|\xi|^\ell
\left(
e^{-t|\xi|^2}
\hat{v}
-\sum_{|\alpha|\le[\gamma]}
M_\alpha(v)
(i\xi)^\alpha e^{-t|\xi|^2}
\right)
\right\|_2
=0,
\end{align}
which has already been stated in Lemma~\ref{5.3}. Furthermore, we can also prove Theorem~\ref{thm:1} by using \eqref{7.3} directly after obtaining \eqref{5.5}. 
As seen in \cite{IKM}, 
sophisticated techniques are needed to prove \eqref{7.1} and \eqref{7.2}. 
In this paper, however, we need only estimates for the $L^2$ framework and so there is room for giving simpler proofs.   
It should be emphasized that independently from the expanding theory established in \cite{IKM}, 
one can actually give alternative proofs based on Lemma~\ref{5.1} which is a generalized result in \cite{I-0}, 
and Lemma~\ref{5.3}. 
This is one of our novelties.  

\vspace{0.3cm}
\noindent{\em Acknowledgement.}
\smallskip
The work of the first author (R. Ikehata) was supported in part by Grant-in-Aid for Scientific Research (C) 15K04958 of JSPS.

\bibliographystyle{amsplain}

\end{document}